\pdfoutput=1
\documentclass[reqno]{amsart}

\newtheorem{theorem}{Theorem}[section]
\newtheorem{lemma}[theorem]{Lemma}
\newtheorem{e-proposition}[theorem]{Proposition}

\newtheorem{e-definition}[theorem]{Definition}

\usepackage{amssymb}
\usepackage{amsmath}
\usepackage{mathrsfs}
\usepackage[all]{xy}

\usepackage{color}
\definecolor{red}{rgb}{1,0,0}

\DeclareMathOperator{\sgn}{sgn}

\DeclareMathOperator{\HKR}{HKR}

\DeclareMathOperator{\End}{End}
\DeclareMathOperator{\complexes}{Ch}
\DeclareMathOperator{\poly}{poly}

\newcommand{\tildeotimes}{\tilde{\otimes}}
\newcommand{\abs}[1]{\left| #1 \right|} 
\newcommand{\into}{\hookrightarrow}

\newcommand{\xto}[1]{\xrightarrow{#1}}

\newcommand{\CC}{\mathbb{C}}

\newcommand{\shuffle}[2]{\mathfrak{S}_{#1}^{#2}}
\newcommand{\sections}[1]{\Gamma(#1)}
\newcommand{\Dpoly}[1]{\mathcal{D}_{\poly}^{#1}}
\newcommand{\enveloppante}[1]{\mathcal{U}(#1)}
\newcommand{\lie}[2]{[#1,#2]} 
\newcommand{\anchor}{\rho}
\newcommand{\category}[1]{\mathcal{#1}}
\newcommand{\groupoid}[1]{\mathscr{#1}}
\newcommand{\bazar}{S^\bullet\big(L/A[-1]\big)}

\begin{document}

\title{A Hopf algebra associated to a Lie pair}

\thanks{{Research partially supported by NSF grant DMS1101827, NSA grant H98230-12-1-0234, 
and NSFC grants 11001146 and 11471179.}}

\author{Zhuo Chen}
\address{Department of Mathematics, Tsinghua University, China}
\email{zchen@math.tsinghua.edu.cn}

\author{Mathieu Sti\'enon}
\address{Department of Mathematics, Penn State University, United States}
\email{stienon@psu.edu}

\author{Ping Xu}
\address{Department of Mathematics, Penn State University, United States}
\email{ping@math.psu.edu}

\maketitle

\begin{center}
\textit{In tribute to Alan Weinstein on the occasion of his seventieth birthday}
\end{center}

\begin{abstract}
The quotient $L/A[-1]$ of a pair $A\into L$ of Lie algebroids 
is a Lie algebra object in the derived category 
$D^b(\category{A})$ of the category $\category{A}$ of left 
$\enveloppante{A}$-modules, the Atiyah class $\alpha_{L/A}$ 
being its Lie bracket. 
In this note, we describe the universal enveloping algebra 
of the Lie algebra object $L/A[-1]$ and we prove that 
it is a Hopf algebra object in $D^b(\category{A})$.
\end{abstract}

\section{Introduction}

Let $A$ be a Lie algebroid over a manifold $M$.
Its space of smooth sections $\sections{A}$ 
is a Lie-Rinehart algebra over the commutative ring $R=C^\infty(M)$.
By an $A$-module, we mean a module over the Lie-Rinehart algebra 
corresponding to the Lie algebroid $A$, i.e. a module over the associative
algebra $\enveloppante{A}$. 

Recall that the universal enveloping algebra $\enveloppante{A}$ 
of a Lie algebroid $A$ over $M$ is simultaneously an associative algebra 
and an $R$-bimodule. 
In case the Lie algebroid $A$ is real, $\enveloppante{A}$ 
is canonically identified to the algebra of left-invariant 
s-fiberwise differential operators on the local Lie groupoid
$\mathscr{A}$ integrating $A$.
Let us recall its construction.

The vector space $\mathfrak{g}=R\oplus\sections{A}$ admits a natural 
Lie algebra structure given by the Lie bracket 
\[ \lie{f+X}{g+Y}=\anchor(X)g-\anchor(Y)f+\lie{X}{Y} ,\] 
where $f,g\in R$ and $X,Y\in\sections{A}$. 
Here $\anchor$ denotes the anchor map. 
Let $i$ denote the natural inclusion of $\mathfrak{g}$ into 
its universal enveloping algebra $\enveloppante{\mathfrak{g}}$. 
The universal enveloping algebra $\enveloppante{A}$ of the Lie algebroid $A$ 
is the quotient of the subalgebra of $\enveloppante{\mathfrak{g}}$ 
generated by $i(\mathfrak{g})$ by the two-sided ideal 
generated by the elements of the form 
$i(f)\otimes i(g+Y)-i(fg+fY)$ with $f,g\in R$ and $Y\in\sections{A}$.

When $A$ is a Lie algebra, $\enveloppante{A}$ is indeed the usual
universal enveloping algebra. On the other hand, when $A$ is the tangent
bundle $TM$, $\enveloppante{A}$ is the algebra of differential operators on $M$.

We use the symbol $\category{A}$ to denote the abelian category 
of $A$-modules. Abusing terminology, we say that a vector bundle 
$E$ over $M$ is an $A$-module if $\sections{E}\in\category{A}$.

Given a Lie pair $(L,A)$ of algebroids, \emph{i.e.} a Lie algebroid $L$ 
with a Lie subalgebroid $A$, the Atiyah class $\alpha_E$ of an $A$-module $E$ 
relative to the pair $(L,A)$ is defined as the obstruction to the existence of 
an \emph{$A$-compatible} $L$-connection on the vector bundle $E$. 
An $L$-connection $\nabla$ on an $A$-module $E$ is said to be 
$A$-compatible if it extends the given flat $A$-connection on $E$ 
and satisfies $\nabla_a\nabla_l-\nabla_l\nabla_a=\nabla_{[a,l]}$ 
for all $a\in\sections{A}$ and $l\in\sections{L}$. 
This fairly recently defined class (see~\cite{CSX}) has as double origin, 
which it generalizes, the Atiyah class of holomorphic vector bundles 
and the Molino class of foliations. 

The quotient $L/A$ of any Lie pair $(L,A)$ is an $A$-module~\cite{CSX}. 
Its Atiyah class $\alpha_{L/A}$ can be described as follows. 
Choose an $L$-connection $\nabla$ on $L/A$ extending the $A$-action. 
Its curvature is the vector bundle map $R^\nabla:\wedge^2 L\to\End(E)$
defined by $R^\nabla(l_1,l_2)=\nabla_{l_1}\nabla_{l_2}-\nabla_{l_2}\nabla_{l_1}
-\nabla_{\lie{l_1}{l_2}}$, for all $l_1, l_2\in\sections{L}$.
Since $L/A$ is an $A$-module, $R^\nabla$ vanishes on $\wedge^2 A$ and, 
therefore, determines a section $R^\nabla_{L/A}$ of 
$A^*\otimes(L/A)^*\otimes\End(L/A)$. 
It was proved in~\cite{CSX} that $R^\nabla_{L/A}$ is a $1$-cocycle 
for the Lie algebroid $A$ with values in the $A$-module 
$(L/A)^*\otimes\End(L/A)$ and that its cohomology class 
$\alpha_{L/A}\in H^1\big(A;(L/A)^*\otimes\End(L/A)\big)$ 
is independent of the choice of the connection. 

Let $\complexes^b(\category{A})$ denote the category of bounded  
complexes in $\category{A}$ and let $D^b(\category{A})$ denote 
the corresponding derived category. 
We write $L/A[-1]$ to denote the quotient $L/A$ regarded 
as a complex in $\category{A}$ concentrated in degree 1. 

The following was proved in~\cite{CSX}.

\begin{e-proposition}[\cite{CSX}]
Let $(L,A)$ be a Lie algebroid pair. 
The Atiyah class $\alpha_{L/A}$ of the quotient $L/A$ relative to the pair $(L,A)$ 
determines a morphism \[ L/A[-1]\otimes L/A[-1]\to L/A[-1] \]
in the derived category $D^b(\category{A})$ making $L/A[-1]$ 
a Lie algebra object in $D^b(\category{A})$. 
\end{e-proposition}

It is well known that every ordinary Lie algebra $\mathfrak{g}$ 
admits a universal enveloping algebra 
$\enveloppante{\mathfrak{g}}$, which is a Hopf algebra. 
We are thus led to the following natural questions: does there exist 
a universal enveloping algebra for $L/A[-1]$ in $D^b(\category{A})$
and, if so, is it a Hopf algebra object?

In this Note, we give a positive answer to the questions above.
For a complex manifold $X$, the Atiyah class of the Lie pair 
$(T_X\otimes\CC,T_X^{0, 1})$ is simply the usual Atiayh class 
of the holomorphic tangent bundle $T_X$  
recently exploited by Kapranov~\cite{Kapranov}. It was proved
that the universal enveloping algebra of the Lie algebra object
$T_X[-1]$  in $D^b (X)$ is the Hochschild cochain complex 
$(\Dpoly{\bullet}(X), d)$~\cite{Markarian,Ramadoss,RW}.
This result played an important role in the study of several aspects of 
complex geometry including the Riemann-Roch theorem~\cite{Markarian},
the Chern character~\cite{Ramadoss} and the Rozansky-Witten 
invariants~\cite{RW,RozanskyWitten}. Applications of our result
will be developed elsewhere.

\section{Hochschild-Kostant-Rosenberg map}

It is known~\cite{Xu:quantum} 
that the universal enveloping algebra $\enveloppante{L}$ 
of a Lie algebroid $L$ admits a cocommutative coassociative coproduct 
$\Delta:\enveloppante{L}\to\enveloppante{L}\tildeotimes\enveloppante{L}$, 
which is defined on generators as follows: 
$\Delta(f)=f\tildeotimes 1=1\tildeotimes f$, $\forall f\in R$ 
and $\Delta(l)=l\tildeotimes 1+1\tildeotimes l$, $\forall l\in\sections{L}$. 
Here, and in the sequel, $\tildeotimes$ stands for the tensor product 
of left $R$-modules. Moreover, $\enveloppante{L}$ is an $L$-module 
since each section $l$ of $L$ acts on $\enveloppante{L}$ 
by left multiplication: $\nabla_l u=l\cdot u$, $\forall u\in\enveloppante{L}$.

Now, given a Lie pair $(L,A)$, consider the quotient $\Dpoly{1}$
of $\enveloppante{L}$ by the left ideal generated by $\sections{A}$. 
It is straighforward to see that the comultiplication on $\mathcal{U}(L)$ 
induces a comultiplication
$\Delta:\Dpoly{1}\to\Dpoly{1}\tildeotimes\Dpoly{1}$ on $\Dpoly{1}$
and the action of $L$ on $\enveloppante{L}$ determines an action of
$A$ on $\Dpoly{1}$.

\begin{lemma}
The quotient 
$\Dpoly{1}=\frac{\enveloppante{L}}{\enveloppante{L}\sections{A}}$ 
is simultaneously a cocommutative coassociative $R$-coalgebra 
and an $A$-module. 
Moreover, its comultiplication is compatible with its $A$-action: 
\[ \nabla_X(\Delta p)=\Delta(\nabla_X p), \quad\forall X\in\sections{A},
p\in\Dpoly{1}. \]
\end{lemma}

Let $\Dpoly{n}$ denote the $n$-th tensorial power 
$\Dpoly{1}\tildeotimes\cdots\tildeotimes\Dpoly{1}$ of $\Dpoly{1}$ 
and, for $n=0$, set $\Dpoly{0}=R$. We define a coboundary operator 
$d:\Dpoly{\bullet}\to\Dpoly{\bullet+1}$ 
on $\Dpoly{\bullet}=\bigoplus_{n=0}^{\infty}\Dpoly{n}$ by 
\begin{multline}\label{eq:d}
d(p_1\tildeotimes\cdots\tildeotimes p_n) 
= 1\tildeotimes p_1\tildeotimes\cdots\tildeotimes p_n 
-(\Delta p_1)\tildeotimes\cdots\tildeotimes p_n
+p_1\tildeotimes(\Delta p_2)\tildeotimes\cdots\tildeotimes p_n -\cdots \\ 
+(-1)^np_1\tildeotimes\cdots\tildeotimes p_{n-1}\tildeotimes(\Delta p_n) 
+(-1)^{n+1}p_1\tildeotimes\cdots\tildeotimes p_n\tildeotimes 1
,\end{multline}
for any $p_1,p_2,\dots,p_n\in\Dpoly{1}$. 
Since the comultiplication $\Delta$ is compatible with the action of $A$, 
the operator $d$ is a morphism of $A$-modules.  
Moreover, $\Delta$ being coassociative, $d$ satisfies $d^2=0$. 
Thus $(\Dpoly{\bullet},d)$ is an object of $\complexes^b(\category{A})$.

When endowed with the trivial coboundary operator, 
the space of sections of 
\[ S^\bullet\big(L/A[-1]\big)=\bigoplus_{k=0}^{\infty} S^k\big(L/A[-1]\big) 
=\bigoplus_{k=0}^{\infty} \big(\wedge^k L/A\big)[-k] \]
is a complex of $A$-modules: 
\[ 0\to R\xto{0} \sections{L/A}\xto{0} \sections{\wedge^2(L/A)}\xto{0} \sections{\wedge^3(L/A)}\xto{0} \cdots \] 
The natural inclusion $\sections{L/A}\hookrightarrow\Dpoly{1}$
extends naturally to the Hochschild-Kostant-Rosenberg map
\[ \HKR:\sections{\bazar}\to\Dpoly{\bullet} \] by skew-symmetrization: 
\begin{multline}\label{eq:hkr} 
\HKR(b_1\wedge\cdots\wedge b_n)=\frac{1}{n!}\sum_{\sigma\in S_n}\sgn(\sigma) 
b_{\sigma(1)}\tildeotimes b_{\sigma(2)}\tildeotimes\cdots\tildeotimes b_{\sigma(n)},
\\ \forall b_1,\cdots,b_n\in\sections{L/A}
.\end{multline}

\begin{e-proposition}\label{Prop:HKRisQIS}
In $\complexes^b(\category{A})$, the Hochschild-Kostant-Rosenberg map is a quasi-isomorphism 
from $(\sections{\bazar},0)$ to $(\Dpoly{\bullet},d)$.
\end{e-proposition}

\paragraph*{Sketch of proof}
Assuming $L$ and $A$ are real Lie algebroids, 
let $\groupoid{L}$ and $\groupoid{A}$ be local Lie groupoids integrating 
$L$ and $A$ respectively. The source map $s:\groupoid{L}\to M$ induces 
a surjective submersion $J:\groupoid{L}/\groupoid{A}\to M$. 
The right quotient $\groupoid{L}/\groupoid{A}$ is a left 
$\groupoid{L}$-homogeneous space with momentum map $J$~\cite{LWX:CMP}. 
Therefore, it admits an infinitesimal $L$-action, and hence an infinitesimal 
$A$-action. The coalgebra $\Dpoly{1}$ may be regarded as the space 
of distributions on the $J$-fibers of $\groupoid{L}/\groupoid{A}$ supported 
on $M$. Its $A$-module structure then stems from the infinitesimal $A$-action 
on $\groupoid{L}/\groupoid{A}$. The $n$-th tensorial power $\Dpoly{n}$ 
may be viewed as the space of $n$-differential operators on the $J$-fibers 
of $\groupoid{L}/\groupoid{A}$ evaluated along $M$ and the differential $d$ 
as the Hochschild coboundary. The conclusion follows from the classical 
Hochschild-Kostant-Rosenberg theorem. To prove the proposition 
for complex Lie algebroids, it suffices to consider formal groupoids 
instead of local Lie groupoids~\cite{KP}. 

\section{Universal enveloping algebra of $L/A[-1]$ in $D^b(\category{A})$}

Following Markarian~\cite{Markarian}, Ramadoss~\cite{Ramadoss}, 
and Roberts-Willerton~\cite{RW}, we introduce the following:
\begin{e-definition}
If it exists, the universal enveloping algebra of a Lie algebra object 
$\mathcal{G}$ in $D^b(\category{A})$ is an associative algebra object 
$\mathcal{H}$ in $D^b(\category{A})$ together with a morphism of Lie algebras 
$i:\mathcal{G}\to\mathcal{H}$ satisfying the following universal property: 
given any associative algebra object $\mathcal{K}$ and any morphism of 
Lie algebras $f:\mathcal{G}\to\mathcal{K}$ in $D^b(\category{A})$, 
there exists a unique morphism of associative algebras 
$f':\mathcal{H}\to\mathcal{K}$ in $D^b(\category{A})$ such that $f=f'\circ i$. 
\end{e-definition}

In view of the similarity between $(\Dpoly{\bullet},d)$ and the Hochschild cochain complex, 
we define a cup product $\cup$ on $\Dpoly{\bullet}$ by setting 
$P\cup Q=P\tildeotimes Q$, for all $P, Q\in\Dpoly{\bullet}$. 
Is is simple to check that
\[ d(P\cup Q) = d P \cup Q + (-1)^{\abs{P}}P\cup d Q ,\] 
for all homogeneous $P,Q\in\Dpoly{\bullet}$.

\begin{e-proposition}\label{pro:2.1}
For any Lie pair $(L,A)$ of algebroids, 
$(\Dpoly{\bullet},d,\cup)$ is an associative algebra
object in $D^b(\category{A})$,
which is in fact the universal enveloping algebra of
the Lie algebra $L/A[-1]$ in $D^b(\category{A})$.
\end{e-proposition}

Consider the inclusion $\eta:R\hookrightarrow\Dpoly{n}$,
the projection $\varepsilon:\Dpoly{n}\twoheadrightarrow R$, 
and the maps $t:\Dpoly{\bullet}\to\Dpoly{\bullet}$ and 
$\tilde{\Delta}:\Dpoly{\bullet}\to\Dpoly{\bullet}\underset{R}{\otimes}\Dpoly{\bullet}$ 
defined, respectively, by 
\[ t(p_1\tildeotimes p_2\tildeotimes\cdots\tildeotimes p_n)
=(-1)^{\frac{n(n-1)}{2}} p_n\tildeotimes p_{n-1}\tildeotimes\cdots\tildeotimes p_1 \]
and 
\[ \tilde{\Delta}(p_1\tildeotimes p_2\tildeotimes\cdots\tildeotimes p_n)= 
\sum_{i+j=n}\sum_{\sigma\in\shuffle{i}{j}}\sgn(\sigma) 
\big(p_{\sigma(1)}\tildeotimes\cdots\tildeotimes p_{\sigma(i)}\big)
\otimes\big(p_{\sigma(i+1)}\tildeotimes\cdots\tildeotimes p_{\sigma(n)}\big) ,\]
where $\shuffle{i}{j}$ denotes the set of
$(i,j)$-shuffles.\footnote{An $(i,j)$-shuffle is a permutation
$\sigma$ of the set $\{1,2,\cdots,i+j\}$ such that
$\sigma(1)\le\sigma(2)\le\cdots\le\sigma(i)$ and
$\sigma(i+1)\le\sigma(i+2)\le\cdots\le\sigma(i+j)$.}

\begin{theorem}\label{thm:2.2}
For any Lie pair $(L,A)$ of algebroids, 
$(\Dpoly{\bullet}, d)$  with the multiplication $\cup$,
the comultiplication $\tilde{\Delta}$, the unit $\eta$, 
the counit $\varepsilon$, and the antipode $t$, 
is a Hopf algebra object in $D^b(\category{A})$. 
\end{theorem}

\section{Ramadoss's approach: $L(\Dpoly{1})$}

To prove Proposition~\ref{pro:2.1} and Theorem~\ref{thm:2.2},
we essentially follow Ramadoss's approach~\cite{Ramadoss}.
Let $L(\Dpoly{1})$ be the (graded) free Lie algebra generated over $R$ 
by $\Dpoly{1}$ concentrated in degree $1$. 
In other words, $L(\Dpoly{1})$ is the smallest Lie subalgebra of $\Dpoly{\bullet}$ 
containing $\Dpoly{1}$. The Lie bracket of two vectors
$u\in\Dpoly{i}$ and $v\in\Dpoly{j}$ is the vector
$[u,v]=u\tildeotimes v-(-1)^{ij}v\tildeotimes u\in\Dpoly{i+j}$.
Actually, $L(\Dpoly{1})$ is made of all linear combinations of 
elements of the form $[p_1,[p_2,[\cdots,[p_{n-1},p_n]\cdots]]]$ with 
$p_1,\dots,p_n\in\Dpoly{1}$. 
One checks that $L(\Dpoly{1})$ is a $d$-stable $A$-submodule of $\Dpoly{\bullet}$ 
and that its Lie bracket is a chain map with respect to the coboundary operator $d$. 
Therefore $(L(\Dpoly{1}),d)$ is a Lie algebra object in $\complexes^b(\category{A})$.

Let $S^{\bullet}(L(\Dpoly{1}))$ be the symmetric algebra of $L(\Dpoly{1})$ 
and let \[ I:S^{\bullet}(L(\Dpoly{1}))\to\Dpoly{\bullet} \] be the symmetrization map:
\[ I(z_1\odot\cdots\odot z_n)
=\frac{1}{n!}\sum_{\sigma\in S_n}\sgn(\sigma;z_1,\cdots,z_n) 
z_{\sigma(1)}\tildeotimes z_{\sigma(2)}\tildeotimes\cdots\tildeotimes z_{\sigma(n)} .\]
The Koszul sign $\sgn(\sigma;z_1,\cdots,z_n)$ of a permutation $\sigma$ of
the (homogeneous) vectors $z_1,z_2,\dots,z_n\in S^{\bullet}(L(\Dpoly{1}))$ 
is determined by the relation
\[ z_{\sigma(1)}\odot z_{\sigma(2)}\odot\cdots\odot z_{\sigma(n)} =
\sgn(\sigma;z_1,\cdots,z_n)  z_1\odot z_2\odot\cdots\odot z_n .\]

\begin{lemma}\label{lem:1}
The symmetrization $I:S^{\bullet}(L(\Dpoly{1}))\to\Dpoly{\bullet}$ 
is an isomorphism in $\complexes^b(\category{A})$.
\end{lemma}

Using Lemma~\ref{lem:1} and the HKR quasi-isomorphism, 
one can prove that the composition $\beta:\sections{L/A[-1]}\to L(\Dpoly{1})$
of the inclusions \[ \sections{L/A[-1]}\subset\Dpoly{1}\subset L(\Dpoly{1}) \] 
is a quasi-isomorphism in $\complexes^b(\category{A})$, 
which intertwines the Lie brackets on $\sections{L/A[-1]}$ and $L(\Dpoly{1})$.

\begin{e-proposition}
\begin{enumerate}
\item The inclusion $\beta:\sections{L/A[-1]}\to L(\Dpoly{1})$ is a quasi-isomorphism in 
$\complexes^b(\category{A})$.
\item The inclusion $\beta:\sections{L/A[-1]}\to L(\Dpoly{1})$ 
is an isomorphism of Lie algebra objects in $D^b(\category{A})$ 
as the diagram 
\[ \xymatrix{
\sections{L/A[-1]} \tildeotimes \sections{L/A[-1]} \ar[d]_{\alpha_{L/A}} \ar[rr]^{\beta\otimes\beta}
&& L(\Dpoly{1}) \tildeotimes L(\Dpoly{1}) \ar[d]^{[,]} \\
\sections{L/A[-1]} \ar[rr]_{\beta} && L(\Dpoly{1})
} \]
commutes in $D^b(\category{A})$.
\end{enumerate}
\end{e-proposition}

Proposition~\ref{pro:2.1} and Theorem~\ref{thm:2.2} now follow immediately. 

\bibliographystyle{amsplain}
\bibliography{mathscinet}

\end{document}